\documentclass[11pt]{article}
\usepackage{eucal}
\usepackage{amscd}
\usepackage{amssymb}
\usepackage{amsmath}
\usepackage{amsfonts}
\usepackage{amsopn}
\begin{document}
\begin{center}
{\bf Selmer groups of elliptic curves with complex multiplication}\\
\vspace*{5mm}
{\sc By A. SAIKIA}\\
{\it McGill University}\\
{\it Montreal, H3A 2K6 Canada.}\\
 \vspace*{5mm}
\end{center}
\vspace*{5mm} {\small{\it Abstract.}  {Suppose $K$ is an imaginary
quadratic field and $E$ is an elliptic curve over a number field
$F$ with complex multiplication by the ring of integers in $K$.
Let $p$ be a rational prime that splits as
$\mathfrak{p}_{1}\mathfrak{p}_{2}$ in $K$. Let $E_{p^{n}}$ denote
the $p^{n}$-division points on $E$. Assume that $F(E_{p^{n}})$ is
abelian over $K$ for all $n\geq 0$. This paper proves that the
Pontrjagin dual of the $\mathfrak{p}_{1}^{\infty}$-Selmer group of
$E$ over $F(E_{p^{\infty}})$ is a finitely generated free
$\Lambda$-module, where $\Lambda$ is the Iwasawa algebra of
Gal$\big(F(E_{p^{\infty}})/F(E_{\mathfrak{p}_{1}^{\infty}\mathfrak{p}_{2}})\big)$.
It also gives a simple formula for the rank of the Pontrjagin dual
as a $\Lambda$-module.}}\\

\noindent {\bf{Acknowledgment}}. The author is indebted to J.H.
Coates for many helpful suggestions at various stages of this
paper. This paper would not have been possible without his
guidance. The author also thanks S. Howson, L. Fu and the referee
for their comments. The author was supported by Hodge Fellowship
at IHES and CRM/CICMA post doctoral fellowship at McGill
University during the progress of this work.

\begin{center} 1. \underline{\it \large Introduction}\end{center}
Let $K$ be an imaginary quadratic field. Suppose $E$ is an
elliptic curve over a number field $F$ with complex multiplication
by the ring of integers ${\cal{O}}$ in $K$. Let $p\neq 2,3$ denote
a rational prime such that $p{\cal{O}}=
\mathfrak{p}_{1}\mathfrak{p}_{2}$ and assume that $E$ has good
reduction over both $\mathfrak{p}_{1}$ and $\mathfrak{p}_{2}$.
Pick any element $\pi$ of ${\cal{O}}$ such that $\pi {\cal{O}} =
\mathfrak{p}_{1}^{h}$ for some $h\geq 1$. Clearly, there is also
an element $\bar{\pi}$ in ${\cal{O}}$ such that
$\bar{\pi}{\cal{O}} = \mathfrak{p}_{2}^{h}$. Let $L$ be an
algebraic extension of $F$. For $n \geq 0$, the $\pi^{n}$-Selmer
group of $E$ over $L$ is defined as
\begin{equation*}
{Sel}_{\pi^{n}}(E/L) = \mbox{Ker}\Big( H^{1}(L,E_{\pi^{n}})
\longrightarrow \prod\limits_{v}H^{1}(L_{v},E)_{\pi^{n}}\Big),
\end{equation*}
where $v$ runs over all the places of $L$. The
$\mathfrak{p}_{1}^{\infty}$-Selmer group of $E/L$ is defined as
\[{Sel}_{\mathfrak{p}_{1}^{\infty}}(E/L)= \lim_{\substack{\longrightarrow\\ n}}
{Sel}_{\pi^{n}}(E/L),\] where the limit is with respect to the
homomorphisms induced by the natural inclusion of $E_{{\pi}^{n}}$
into $E_{{\pi}^{n+1}}$. The $\mathfrak{p}_{1}^{\infty}$-Selmer
group fits into an exact sequence
\begin{equation} \label{selm} 0 \longrightarrow
E(L) \otimes
K_{\mathfrak{p}_{1}}/{{\cal{O}}_{\mathfrak{p}_{1}}}\longrightarrow
{Sel}_{\mathfrak{p}_{1}^{\infty}}(E/L) \longrightarrow III
(E/L)_{\mathfrak{p}_{1}^{\infty}}\longrightarrow 0,
\end{equation} where $E(L)$ is the Mordell-Weil group of rational
points on $E$ defined over $L$ and $III(E/L)$ is the
Tate-Shafarevich group of $E/L$ defined by
\[III(E/L) = \mbox{Ker} \Big( H^{1}(L,E ) \longrightarrow \prod\limits_{v}
H^{1}(L_{v},E)\Big).\] One of the basic questions in number theory
is to understand the Mordell-Weil group and the Tate-Shafarevich
group of $E$ over various field extensions of $\mathbb{Q}$. Thus,
the importance of the study of Selmer groups arise from the exact
sequence (\ref{selm}) above.\\

There are some natural choices for the field extension $L$ of $F$,
over which we want to examine the structure of
${Sel}_{\mathfrak{p}_{1}^{\infty}}(E/L)$. We usually take $L$ to
be a field generated over $F$ by the torsion points on $E$. In
particular, we will consider
\[F_{\infty}=F(E_{p^{\infty}}),\]
and study ${Sel}_{\mathfrak{p}_{1}^{\infty}}(E/F_{\infty})$, or rather its
Pontrjagin dual $X(F_{\infty})$. By definition,
\[X(F_{\infty})= Hom \big({Sel}_{\mathfrak{p}_{1}^{\infty}}
(E/F_{\infty}), \mathbb{Q}_{p}/\mathbb{Z}_{p}\big).\]
It is compact and has the natural structure of
Gal($F_{\infty}/F$)-module.
This will be the primary object of our study in this paper.\\

\begin{center}2. \underline{\it \large Notation} \end{center}

We define the following field extensions of the number field $F$
generated by torsion points on $E$:
\[L_{0} = F(E_{p}),\;\; F_{0}= L_{0}(E_{\mathfrak{p}_{1}^{\infty}}),\;\;
L_{\infty} = L_{0}(E_{\mathfrak{p}_{2}^{\infty}}),\;\;
F_{\infty}=F(E_{p^{\infty}}).\] Let $\Gamma'$ be the Galois group
of $F_{\infty}$ over $L_{0}$, and $\Sigma$ be the Galois group
$F_{0}$ over $L_{0}$. Let $\Gamma$ be the Galois group
$F_{\infty}$ over $F_{0}$, which can also be identified with the
Galois group $L_{\infty}$ over $L_{0}$. Clearly, $\Gamma'$ is
isomorphic to $\mathbb{Z}_{p}^{2}$, whereas $\Gamma$ and $\Sigma$
are isomorphic to $\mathbb{Z}_{p}$. We denote the unique subgroup
of index $p^{n}$ in $\Gamma$ by $\Gamma_{n}$. Let $L_{n}$ and
$F_{n}$ be the fixed fields of $L_{\infty}$ and $F_{\infty}$
respectively under the action of $\Gamma_{n}$. Then, we have the
following Galois groups:
\[{{\mbox{Gal}}(L_{\infty}/L_{n})\simeq
\mbox{{Gal}}(F_{\infty}/F_{n})= \Gamma_{n}, \;\;\;\;\;
\mbox{{Gal}}(L_{n}/L_{0})\simeq\mbox{{Gal}}
(F_{n}/F_{0})=\Gamma/\Gamma_{n} \simeq
\mathbb{Z}_{p}/p^{n}\mathbb{Z}_{p}}.\] We have the following field
diagram:
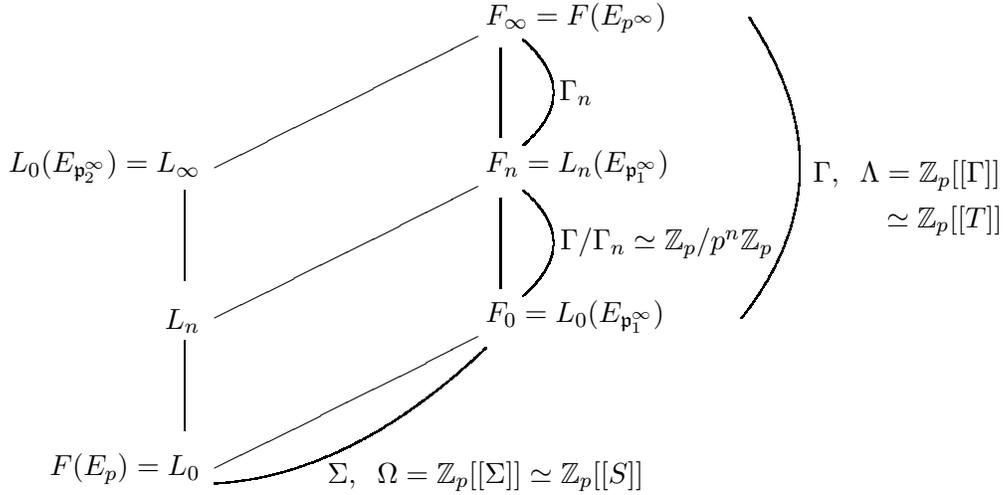
\begin{center}
\setlength{\unitlength}{1mm}
\begin{picture}(110,70)
\put(10,43){\makebox(0,0)[rb]{$L_{0}(E_{\mathfrak{p}_{2}^{\infty}})=L_{\infty}$}}
\put(10,23){\makebox(0,0)[rb]{$L_{n}$}}
\put(10,3){\makebox(0,0)[rb]{$F(E_{p})=L_{0}$}}
\put(48,65){\makebox(0,0)[l]{$F_{\infty}=F(E_{p^{\infty}})$}}
\put(48,43){\makebox(0,0)[lb]{$F_{n}=L_{n}(E_{\mathfrak{p}_{1}^{\infty}})$}}
\put(48,23){\makebox(0,0)[lb]{$F_{0}=L_{0}(E_{\mathfrak{p}_{1}^{\infty}})$}}
\put(12,5){\line(2,1){35}} \put(12,25){\line(2,1){35}}
\put(12,45){\line(2,1){35}} \put(8,10){\line(0,1){12}}
\put(8,30){\line(0,1){12}} \put(50,29){\line(0,1){12}}
\put(50,49){\line(0,1){12}} \qbezier(83,65)(97,44)(82,25)
\qbezier(53,62)(61,55)(53,48) \qbezier(53,42)(61,35)(53,28)
\put(72,35){\makebox(0,0){$\Gamma/\Gamma_{n}\simeq
\mathbb{Z}_{p}/p^{n}\mathbb{Z}_{p}$}}
\put(104,44){\makebox(0,0){$\Gamma,\;\;\Lambda =
\mathbb{Z}_{p}[[\Gamma]]$}} \put(109,38){\makebox(0,0){$\simeq
\mathbb{Z}_{p}[[T]]$}} \put(60,55){\makebox(0,0){$\Gamma_{n}$}}
\qbezier(12,3)(30,4)(48,21) \put(48,4){\makebox(0,0){$\Sigma,\;\;
\Omega = \mathbb{Z}_{p}[[\Sigma]]\simeq \mathbb{Z}_{p}[[S]]$}}
\end{picture}
\end{center}
The Iwasawa algebra of $\Gamma$ is defined as
\[\mathbb{Z}_{p}[[\Gamma]] = \lim_{\substack{\leftarrow \\ n}}\mathbb{Z}_{p}[\Gamma/\Gamma_{n}],\]
where the inverse limit is taken with respect to canonical
surjective maps. We denote the Iwasawa algebra of $\Gamma$ by
$\Lambda $, and that of $\Sigma$ by $\Omega$. Following Serre, we
can identify $\Lambda$ with $\mathbb{Z}_{p}[[T]]$ and $\Omega$
with $\mathbb{Z}_{p}[[S]]$. We note that
$\mathbb{Z}_{p}[[\Gamma']]$ is isomorphic to
$\mathbb{Z}_{p}[[T,S]]$. We will denote the Pontrjagin dual of
$Sel_{\mathfrak{p}_{1}^{\infty}}(E/F_{n})$ by $X(F_{n})$.\\

\begin{center}3. \underline{\it \large Statement of results} \end{center}

Our goal is to study the structure of $X(F_{\infty})$ as a module
over the Iwasawa algebra $\Lambda \simeq \mathbb{Z}_{p}[[T]]$. We
shall work under
the following hypothesis : \\

\noindent ({\bf Hyp}) The fields $F(E_{p^{n}})$ are abelian over $K$ for all $n\geq 0$. \\

\noindent Note that when $F=K$, the hypothesis is true by theory
of complex multiplication. It is well known (e.g. see [P-R 1])
that $X(F_{\infty})$ is a finitely generated torsion module over
the Iwasawa algebra $\mathbb{Z}_{p}[[S,T]]$, whereas $X(F_{n})$ is
a finitely generated torsion $\mathbb{Z}_{p}[[S]]$-module under
the above hypothesis. Let $\lambda_{0}$ be the rank of $X(F_{0})$
as a $\mathbb{Z}_{p}$-module. In this paper, we shall prove the
following two
theorems about the $\Lambda$-module structure of $X(F_{\infty})$: \\

\noindent {\bf Theorem 1 :} $X(F_{\infty})$ {\it is a finitely generated}
$\Lambda$-{\it module}. \\

\noindent {\bf Theorem 2 :} $X(F_{\infty})$ {\it is a free}
$\Lambda$-{\it module of rank} $\lambda_{0}+r-1$.\\

\noindent Here $r$ is the number of primes of $F_{0}$ above
$\mathfrak{p}_{2}$. Since the primes over $\mathfrak{p}_{2}$ do
not split in the tower $F_{\infty}$ over $F_{0}$, $r$ is also the
number of primes above $\mathfrak{p}_{2}$ of $F_{n}$ for any $n
\geq 0$.\\

\begin{center} 4. \underline{{\it \large The structure of} $X(F_{n})$}
\end{center}

The key idea in the proof of the theorems 1 and 2 is to examine
the relation between $X(F_{\infty})$ and $X(F_{n})$, and then
exploit well-known facts about $X(F_{n})$. Theorem 18 and
proposition 20 in [P-R 1] show that $X(F_{n})$ is a finitely
generated torsion $\mathbb{Z}_{p}[[S]]$-module provided Leopoldt's
conjecture is true for the $\mathbb{Z}_{p}$-extension $F_{n}$ over
$L_{n}$. Brumer proved that Leoplodt's conjecture is true for the
$\mathbb{Z}_{p}$-extensions of an abelian extension of an
imaginary quadratic field. Under our hypothesis (Hyp), $L_{n}$ is
an abelian extension of the imaginary quadratic field $K$.
Therefore, Leopoldt's conjecture holds for $F_{n}$ and as a
consequence, we know that $X(F_{n})$ is a finitely generated
torsion $\mathbb{Z}_{p}[[S]]$-module. By structure theory of
finitely generated torsion $\mathbb{Z}_{p}[[S]]$-module, there is
a homomorphism
\begin{equation}
\label{structure} \phi : X(F_{n}) \longrightarrow  \bigoplus
\Big(\oplus_{i=1}^{s} \mathbb{Z}_{p}[[S]]/p^{n_{i}}\Big) \bigoplus
\Big(\oplus_{j=1}^{t} \mathbb{Z}_{p}[[S]]/(f_{j}^{m_{j}})\Big),
\end{equation}
with finite kernel and cokernel. {Here} $f_{j}$ {are distinguished
polynomials in} $\mathbb{Z}_{p}[[S]]$ {and}
$s,\;t,\;n_{i},\;m_{j}$ {are non-negative integers}. The {\it
lambda}-invariant $\lambda_{n}$ and the {\it mu}-invariant
$\mu_{n}$ of the $\mathbb{Z}_{p}[[S]]$-module $X(F_{n})$ are
defined as
\[\lambda_{n}= \sum_{j=1}^{t}m_{j}.\mbox{deg}(f_{j}),\;\;
\mu_{n} = \sum_{i=1}^{r}n_{i}.\] When $L_{n}$ is an abelian
extension of $K$, Gillard ([Gi 1], [Gi 2]) has shown that $\mu_{n}
=0$. While [Gi 2] has the proof of vanishing of the {\it
mu}-invariant without any assumption on the class number of $K$,
the proof in [Gi 1] works under the assumption that the class
number of $K$ is 1 (that would have amounted to assuming that $E$
is defined over $K$ in our work). As $L_{n}$ is abelian over $K$
under our hypothesis (Hyp), Gillard's result implies that the
$p$-torsion part in the right hand side of (\ref{structure}) does
not occur. Moreover, it follows (as pointed out in theorem 25 of
[P-R 1]) from the work of Greenberg ([Gr 1]) that $X(F_{n})$ has
no finite non-zero $\mathbb{Z}_{p}[[S]]$-submodule. Thus, the
kernel of $\phi$ ({\it a priori} finite) is trivial. Hence, $\phi$
maps $X(F_{n})$ injectively into a free $\mathbb{Z}_{p}$-module of
rank $\lambda_{n}$ with finite cokernel. We have now obtained the
following information regarding the $\mathbb{Z}_{p}$-module
structure
of $X(F_{n})$:\\

\noindent{\bf Proposition 3 : } $X(F_{n})$ {\it is a free}
$\mathbb{Z}_{p}$-{\it module of rank} $\lambda_{n}$ {\it under our
hypothesis (Hyp)}.\\

\noindent How the {\it lambda}-invariant $\lambda_{n}$ of
$X(F_{n})$ varies along the tower of fields $F_{n}$ ($n=0,1,2,
\ldots$) will be very important to us. We will study this question
in section 7 (c.f. lemma 11).

\begin{center} 5. \underline{\it \large A crucial proposition} \end{center}

Let us fix an $n \geq 0$. Let $S$ be the set of primes of
$F$ above $p$. Let $F_{S}$ be the maximal extension of $F$
unramified outside $S$. It is clear that $F_{\infty}\subset F_{S}$
and $E_{\mathfrak{p}_{1}^{\infty}}\subset E(F_{S})$. The following
result is a crucial ingredient in examining the relation between
$X(F_{\infty})$ and $X(F_{n})$ [see the commutative diagram (c.d.)
in section 6] : \\

\noindent{\bf Proposition 4 :} {\it There is an exact sequence of
Galois modules }
\[{0}  \longrightarrow { {Sel}_{\mathfrak{p}_{1}^{\infty}}(E/F_{n}) }
\longrightarrow  {H^{1}(F_{S}/F_{n}, E_{\mathfrak{p}_{1}^{\infty}})}
\longrightarrow { \prod\limits_{v \mid p}
H^{1}(F_{n,v}, E)_{\mathfrak{p}_{1}^{\infty}}}   \longrightarrow  0 .\]

The key part in the above proposition is the surjectivity.
Hachimori and Matsuno [H-M] proved the above result for the
cyclotomic $\mathbb{Z}_{p}$-extension of a number field. But their
argument carries over to our situation of elliptic curves with
complex multiplication. We will briefly describe how the methods
of [H-M] can be adopted in our case. We will see that the sequence
in proposition 4 comes from a five-term Cassels-Poitou-Tate
sequence (\ref{cpt3}). It will be sufficient to show that the
fourth term in (\ref{cpt3}) vanishes (lemma 5). As a consequence
of this method of proof, we deduce that the fifth term in
(\ref{cpt3}) (a $H^{2}$ term) also vanishes and deduce corollary
6. This vanishing (of $H^{2}$) will be needed for the calculations
of section 7, especially lemma 12. \\

Let us denote the $\mathbb{Z}_{p}$-extension $F_{n}$ of $L_{n}$ by
$T_{\infty}$. We know that the Galois group $\Sigma \simeq \mbox
{Gal}(T_{\infty}/L_{n}$) has a unique subgroup $\Sigma_{m}$ of
index $p^{m}$. Let $T_{m}$ be the fixed field of $T_{\infty}$
under the action of $\Sigma_{m}$. We have a field diagram\\

\begin{center}
\setlength{\unitlength}{1mm}
\begin{picture}(50,46)
\put(10,43){\makebox(0,0)[rb]
{$F_{n}=L_{n}(E_{{\mathfrak{{p}}}_{1}^{\infty}})=T_{\infty}$}}
\put(10,23){\makebox(0,0)[rb] {$T_{m}$}}
\put(10,3){\makebox(0,0)[rb]{$L_{n}$}} \put(8,10){\line(0,1){12}}
\put(8,30){\line(0,1){12}} \qbezier(13,41)(18,33)(13,25)
\qbezier(23,5)(32,25)(23,45)
\put(20,36){\makebox(0,0){$\Sigma_{m}$}}
\put(47,25){\makebox(0,0){$\Sigma,\;\;\Omega\simeq
\mathbb{Z}_{p}[[\Sigma]]$}}
\end{picture}
\end{center}

\noindent By Cassels-Poitou-Tate sequence for the number fields
$T_{m}$, we have a long exact sequence (where $\widehat{M}$
denotes the Pontrjagin dual of $M$)
\begin{equation}
\label{cpt1}
\begin{split}
 0 & \longrightarrow {Sel}_{\pi^{k}}(E/T_{m}) \longrightarrow
 {H^{1}(F_{S}/T_{m}, E_{\pi^{k}})}   \longrightarrow
 { \prod\limits_{v \mid p}H^{1}(T_{m,v}, E)_{\pi^{k}}}  \\
& \longrightarrow \widehat{{Sel}_{\bar{\pi}^{k}}(E/T_{m})}
\longrightarrow {H^{2}(F_{S}/T_{m}, E_{\pi^{k}})}
\longrightarrow { \prod\limits_{v \mid p}H^{2}
(T_{m,v}, E_{\bar{\pi}^{k}})}  \\
 & \longrightarrow
 \widehat{H^{0}(F_{S}/T_{m}, E_{\bar{\pi}^{k}})}\longrightarrow 0.
\end{split}
\end{equation}
We note that in applying Poitou-Tate duality, one has to consider
not only the primes above $p$, but also the infinite primes and
the primes of bad reduction. However, $E$ has good reduction
everywhere over $L_{0}$ by theory of complex multiplication, and
we can also ignore the infinite primes as $p$ is odd. The
inclusion $E_{\pi^{k}}\hookrightarrow E_{\pi^{k+1}}$ induces a map
$H^{i}(F_{S}/T_{m}, E_{\pi^{k}})$ to $H^{i}(F_{S}/T_{m},
E_{\pi^{k+1}})$, and its dual is given by `multiplication by
$\pi$'. By taking direct limits in (\ref{cpt1}) as $k$ goes to
infinity, we get a five term exact sequence
\begin{equation}
\label{cpt2}
\begin{split}
0 & \longrightarrow {Sel}_{\mathfrak{p}_{1}^{\infty}}(E/T_{m})
\longrightarrow {H^{1}(F_{S}/T_{m}, E_{\mathfrak{p}_{1}^{\infty}})}
 \longrightarrow { \prod\limits_{v \mid p}
 H^{1}(T_{m,v}, E)_{\mathfrak{p}_{1}^{\infty}}} \\
 & \longrightarrow {\Big({{\lim_{\substack{\longleftarrow\\ k}}
 {{Sel}_{\bar{\pi}^{k}}(E/T_{m})}}}\Big)}^{\wedge}
 \longrightarrow {H^{2}(F_{S}/T_{m}, E_{\mathfrak{p}_{1}^{\infty}})}
 \longrightarrow 0.
\end{split}
\end{equation}
We remark that when we take direct limit with respect to $k$, the
sixth term in (\ref{cpt1}) vanishes by Tate local duality (see
Ch.II prop. 16 in [Se]). There is a restriction map from
$H^{i}(F_{S}/T_{m}, E_{\mathfrak{p}_{1}^{\infty}})$ to
$H^{i}(F_{S}/T_{m+1}, E_{\mathfrak{p}_{1}^{\infty}})$, and the
dual map is given by corestriction which acts like the norm map on
$H^{0}$. We now take direct limits in (\ref{cpt2}) as $m$ goes to
infinity, and obtain a five term exact sequence
\begin{equation}
\label{cpt3}
\begin{split}
0 & \longrightarrow {Sel}_{\mathfrak{p}_{1}^{\infty}}(E/T_{\infty})
\longrightarrow {H^{1}(F_{S}/T_{\infty}, E_{\mathfrak{p}_{1}^{\infty}})}
\longrightarrow { \prod\limits_{v \mid p}
H^{1}(T_{\infty,v}, E)_{\mathfrak{p}_{1}^{\infty}}} \\
 & \longrightarrow {\Big({{\lim_{\substack{\longleftarrow\\ m}}
 \lim_{\substack{\longleftarrow\\ k}}
 {{Sel}_{\bar{\pi}^{k}}(E/T_{m})}}}\Big)}^{\wedge}
 \longrightarrow {H^{2}(F_{S}/T_{\infty}, E_{\mathfrak{p}_{1}^{\infty}})}
 \longrightarrow 0.
\end{split}
\end{equation}
Let us denote the fourth term in the above sequence as $\hat{W}$, i.e.,
\[W =  {{\lim_{\substack{\longleftarrow\\ m}}
\lim_{\substack{\longleftarrow\\ k}}{{Sel}_{\bar{\pi}^{k}}(E/T_{m})}}}.\]

Proposition 4 claims that the fourth term in the above
sequence (\ref{cpt3}) vanishes.\\

\noindent {\bf Lemma 5 :} $$W = \lim_{\substack{\longleftarrow\\
m}} \lim_{\substack{\longleftarrow\\
k}}{{Sel}_{\bar{\pi}^{k}}(E/T_{m})}=0.$$

\noindent{\bf Proof} : We adopt an argument similar to the one in
proposition 2.3 of [H-M]. We have an exact sequence (see lemma 1.8
in [C-S])
\[0 \longrightarrow E_{\bar{\pi}^{\infty}}(T_{m}) \longrightarrow
\lim_{\substack{\longleftarrow\\ k}}
{{Sel}_{\bar{\pi}^{k}}(E/T_{m})} \longrightarrow
Hom_{\mathbb{Z}_{p}}
\big(\widehat{{{Sel}_{\bar{\pi}^{\infty}}(E/T_{m})}},
\mathbb{Z}_{p}\big)\longrightarrow 0.\] We now take inverse limit
with respect to corestriction maps as $m$ goes to infinity. These
maps act like norm maps on the first term, and it vanishes in the
limit since only finitely many $\bar{\pi}$-torsion points of $E$
are defined over $T_{\infty}$. Thus, we obtain an injection
\[W = \lim_{\substack{\longleftarrow\\ m}}
\lim_{\substack{\longleftarrow\\
k}}{{Sel}_{\bar{\pi}^{k}}(E/T_{m})}
\hookrightarrow \lim_{\substack{\longleftarrow\\
m}}Hom_{\mathbb{Z}_{p}} \big(\widehat{
{{Sel}_{\mathfrak{p}_{2}^{\infty}}(E/T_{m})}},
\mathbb{Z}_{p}\big).\] The kernel of the restriction map
${{Sel}_{\mathfrak{p}_{2}^{\infty}}(E/T_{m})} \longrightarrow
{{Sel}_{\mathfrak{p}_{2}^{\infty}} (E/T_{\infty})}^{\Sigma_{m}}$
is finite and its order is bounded independent of $m$ (this kernel
is contained in $H^{1}(\Sigma_{m},
E_{\mathfrak{p}_{2}^{\infty}}(T_{\infty}))$, and this group is
bounded independent of $m$, as shown in lemma 3.1 of [Gr 2]).
Therefore, we have an injection
\[ \lim_{\substack{\longleftarrow\\ m}}Hom_{\mathbb{Z}_{p}}
\big(\widehat{ {{Sel}_{\mathfrak{p}_{2}^{\infty}}(E/T_{m})}},
\mathbb{Z}_{p}\big) \hookrightarrow \lim_{\substack{\longleftarrow\\
m}}Hom_{\mathbb{Z}_{p}} \big((\widehat{
{{Sel}_{\mathfrak{p}_{2}^{\infty}} (E/T_{\infty})}})_{\Sigma_{m}},
\mathbb{Z}_{p}\big).\] The latter module has the same underlying
set as $Hom_{\Omega}\big(\widehat{
{{Sel}_{\mathfrak{p}_{2}^{\infty}} (E/T_{\infty})}}, \Omega\big)$
(e.g., \S 2, lemma 4(ii) in [P-R 2]).

We again invoke proposition 20 in [P-R 1] which says that
${Sel}_{\mathfrak{p}_{2}^{\infty}}(E/T_{\infty})$ is
$\Omega$-cotorsion provided Leopoldt's conjecture is true for the
$\mathbb{Z}_{p}$-extension $T_{\infty}$ of $L_{n}$. But Leopoldt's
conjecture is true for the $\mathbb{Z}_{p}$-extension $T_{\infty}$
of the abelian [under our hypothesis (Hyp)] extension $L_{n}$ of
the imaginary quadratic field $K$. Therefore,
${Sel}_{\mathfrak{p}_{2}^{\infty}}(E/T_{\infty})$ is
$\Omega$-cotorsion and $Hom_{\Omega}\big(\widehat{
{{Sel}_{\mathfrak{p}_{2}^{\infty}} (E/T_{\infty})}}, \Omega\big) =
0$.

Thus, the compact $\Omega$-module $W$ can be embedded into the null module.
\hspace*{1cm}$\square$\\

With this lemma, the proof of proposition 4 is now complete. The
following corollary to lemma 5 will be a vital step in our proof
of theorem 2 (lemma 12 in section 7).\\

\noindent {\bf Corollary 6 :} For any $n \geq 0$, $
H^{2}\big(F_{S}/F_{n}, E_{\mathfrak{p}_{1}^{\infty}}\big) = 0$.\\

\noindent{\bf Proof :} From the Cassels-Poitou-Tate sequence
(\ref{cpt3}) and lemma 5, it is clear that
$H^{2}\big(F_{S}/T_{\infty},E_{\mathfrak{p}_{1}^{\infty}}\big)=0$.
But $T_{\infty}$ stands for any of the $F_{n}$ for $n \geq 0$.
\hspace*{1cm}$\square$\\

\begin{center}6. \underline{{\it \large Relation between} $X(F_{\infty})$ {\it and}
$X(F_{n})$}
\end{center}

In order to examine the relation between $X(F_{\infty})$ and
$X(F_{n})$, the following commutative diagram is of crucial
importance:

\[\begin{CD}
\scriptscriptstyle {0} @>>> \scriptscriptstyle
{ {Sel}_{\mathfrak{p}_{1}^{\infty}}(E/F_{\infty})^{\Gamma_{n}}}
@>>> \scriptscriptstyle { H^{1}(F_{S}/F_{\infty},
E_{\mathfrak{p}_{1}^{\infty}})^{\Gamma_{n}}}  @>>>
\scriptscriptstyle { \prod\limits_{v \mid p}
\big(\prod\limits_{w\mid v}H^{1}
(F_{\infty,w}, E)_{\mathfrak{p}_{1}^{\infty}}\big)^{\Gamma_{n}}}\\
 @.           @A\scriptscriptstyle {\alpha_{n}}AA
 @A\scriptscriptstyle {\beta_{n}}AA
 @AA\scriptscriptstyle {\gamma_{n}=\prod_{v\mid p} \gamma_{n,v}}A \\
\scriptscriptstyle {0}  @>>>  \scriptscriptstyle {
{Sel}_{\mathfrak{p}_{1}^{\infty}}(E/F_{n}) }  @>>>
\scriptscriptstyle {H^{1}(F_{S}/F_{n},
E_{\mathfrak{p}_{1}^{\infty}})}   @>>>
\scriptscriptstyle { \prod\limits_{v \mid p}
H^{1}(F_{n,v}, E)_{\mathfrak{p}_{1}^{\infty}}}
@>>>\scriptscriptstyle { 0 }\\
\end{CD}\]
\begin{center} Commutative Diagram (c.d.)\end{center}

\noindent  The horizontal maps originate from Cassels-Poitou-Tate
sequence, whereas the vertical maps are induced by restriction.
All of our work in section 5 has been to establish the exactness
of the bottom row in the above diagram. We are primarily
interested in the kernel and cokernel of the map $\alpha_{n}$
above. By the snake lemma, we have an exact sequence
\begin{equation}
\label{snake1} 0 \longrightarrow \mbox{Ker}(\alpha_{n})
\longrightarrow \mbox{Ker}(\beta_{n})\longrightarrow
\mbox{Ker}(\gamma_{n}) \longrightarrow
\mbox{Coker}(\alpha_{n})\longrightarrow
\mbox{Coker}(\beta_{n})\ldots \end{equation} In order to
understand the structure of $\mbox{Ker}(\alpha_{n})$ and
$\mbox{Coker}(\alpha_{n})$, we will first study the kernels and
cokernels of the maps $\beta_{n}$ and $\gamma_{n}$.\\

\noindent  {\bf Lemma 7 :}
{Ker$(\beta_{n})\simeq\mathbb{Q}_{p}/\mathbb{Z}_{p}$, and
Coker$(\beta_{n})=0$}.\\

\noindent {\bf {Proof :}} Recall that all the points in
$E_{\mathfrak{p}_{1}^{\infty}}$ are defined over $F_{n}$ ($n=0,1,
\ldots$). By the inflation-restriction sequence of cohomology,
Ker$(\beta_{n})$ equals $H^{1}(\Gamma_{n},
E_{\mathfrak{p}_{1}^{\infty}})$, and Coker$(\beta_{n})$ is
contained in $H^{2} (\Gamma_{n}, E_{\mathfrak{p}_{1}^{\infty}})$.
But $\Gamma_{n}$ is isomorphic to $\mathbb{Z}_{p}$, and hence it
has $p$-cohomological dimension 1. Therefore, $H^{2} (\Gamma_{n},
E_{\mathfrak{p}_{1}^{\infty}})$ vanishes and it follows that
Coker$(\beta_{n})$ is trivial. Moreover, $\Gamma_{n}$ acts
trivially on $E_{\mathfrak{p}_{1}^{\infty}}$ and hence
$H^{1}(\Gamma_{n}, E_{\mathfrak{p}_{1}^{\infty}})$ equals $Hom
(\Gamma_{n}, \mathbb{Q}_{p}/{\mathbb{Z}}_{p})$. We can now
conclude that Ker$(\beta_{n}) \simeq
\mathbb{Q}_{p}/{\mathbb{Z}}_{p}$.
\hspace*{1cm} $\square$\\

\noindent  {\bf Lemma 8 :} {\it For}
$v|\mathfrak{p}_{1}$, {\it Ker}$(\gamma_{n,v}) =  0 $.\\

\noindent We shall give a short and direct proof of this lemma,
though it follows from a more general result of Perrin-Riou (lemma
9 in [P-R 1]).\\

\noindent {\bf Proof :}  By Shapiro's lemma,
\[\Big(\prod\limits_{w\mid v}H^{1}
(F_{\infty,w}, E)\Big)^{\Gamma_{n}}_{\mathfrak{p}_{1}^{\infty}} =
H^{1}(F_{\infty,w},
E)^{\Gamma_{n,v}}_{\mathfrak{p}_{1}^{\infty}},\] where
$\Gamma_{n,v}$ is the decomposition subgroup of $\Gamma_{n}$. By
the inflation-restriction sequence,
\[\mbox{Ker}(\gamma_{n,v})=H^{1}\big(\Gamma_{n,v},
E(F_{\infty,w})\big)_{\mathfrak{p}_{1}^{\infty}}.\] Clearly,
\[F_{\infty, w} = \bigcup_{{M}} L_{\infty, v'}M, \]
where $M$ runs over the finite extensions of $L_{n,\tilde{v}}$
contained in $F_{n,v}$, and $v',\; \tilde{v}$ are the primes below
$w$ of $L_{\infty}$ and $L_{n}$ respectively. Now,
\[\mbox{Ker}(\gamma_{n,v}) = \lim_{\substack{\rightarrow \\ M}}
H^{1}\big(G(L_{\infty,v'}M/M),
E(L_{\infty,v'}M)\big)_{\mathfrak{p}_{1}^{\infty}}.\] Note that
$E$ has good reduction over $L_{n,\tilde{v}}$. Therefore,
$L_{\infty,v'}$ is unramified over $L_{n,\tilde{v}}$ and so is
$L_{\infty, v'}M $ over $M$. Hence,
${H^{1}\big(G(L_{\infty,v'}M/M), E(L_{\infty,v'}M)\big) = 0}$ (see
[Mi, p. 58]). This concludes the proof of lemma 8. \hspace*{1cm} $\square$  \\

\noindent {\bf Lemma 9 :} {\it For} $v|\mathfrak{p}_{2}$,
{\it Ker}$(\gamma_{n,v}) \simeq \mathbb{Q}_{p}/{\mathbb{Z}}_{p}$.\\

\noindent {\bf Proof :} The extension $F_{\infty}$ is totally
ramified over $F_{n}$ at the prime $v$ over $\mathfrak{p}_{2}$.
Therefore, there is only one prime $w$ of $F_{\infty}$ over $v$
and the decomposition group $\Gamma_{n,v}$ is the Galois group
$\Gamma_{n}$. By the inflation-restriction sequence,
\[\mbox{Ker}(\gamma_{n,v})=H^{1}\big(\Gamma_{n,v},
E(F_{\infty,v})\big)_{\mathfrak{p}_{1}^{\infty}}.\]
 Let $\mathfrak{m}_{\infty,v}$ be the maximal ideal of $F_{\infty,v}$ and
 $k_{\infty,v}$ be the residue field. Let $\hat{E}$ be the formal
 group attached to $E$ giving the kernel of reduction at $v$. We
 have the following exact sequence of $\Gamma_{n,v}$-modules:
\[ 0 \longrightarrow \hat{E}(\mathfrak{m}_{\infty,v}) \longrightarrow
E(F_{\infty , v}) \longrightarrow \tilde{E}_{v}(k_{\infty,v}) \longrightarrow 0. \]
Taking Galois cohomology, we get the following exact sequence:
\begin{equation}
\nonumber
\begin{split}
\ldots \longrightarrow
H^{1}\big(\Gamma_{n,v},\hat{E}(\mathfrak{m}_{\infty,v})\big)_{\mathfrak{p}_{1}^{\infty}}
\longrightarrow
H^{1}\big(\Gamma_{n,v}, E(F_{\infty , v})\big)_{\mathfrak{p}_{1}^{\infty}}\\
\longrightarrow
H^{1}\big(\Gamma_{n,v},\;\tilde{E}_{v}(k_{\infty,v})\big)_{\mathfrak{p}_{1}^{\infty}}
\longrightarrow
H^{2}\big(\Gamma_{n,v},\hat{E}(\mathfrak{m}_{\infty,v})\big)_{\mathfrak{p}_{1}^{\infty}}\rightarrow
\ldots .
\end{split}
\end{equation}
Since $v|\mathfrak{p}_{2}$, $\pi$ is an automorphism of $\hat{E}$.
Therefore,
$H^{i}\big(\Gamma_{n,v},\;\hat{E}(\mathfrak{m}_{\infty,v})\big)_{\mathfrak{p}_{1}^{\infty}}$=$0
\;\forall \; i\geq 0$. Hence we have
\[ H^{1}\big(\Gamma_{n,v},\; E(F_{\infty , v})\big)_{\mathfrak{p}_{1}^{\infty}}
\stackrel{\sim}{\longrightarrow} H^{1}\big(\Gamma_{n,v},\;
\tilde{E}_{v}(k_{\infty,v})\big)_{\mathfrak{p}_{1}^{\infty}}.\] As
$\tilde{E}_{v}(k_{\infty,v})$ is a torsion module, we can take the
$\mathfrak{p}_{1}^{\infty}$-torsion inside the cohomology group.
Since $F_{\infty,w}$ is totally ramified over $F_{n,v}$, the group
$\Gamma_{n,v}$ acts trivially on $\tilde{E}_{v}(k_{\infty,v})$.
Therefore, the right hand side in the previous expression is
\[Hom(\Gamma_{n,v}, \tilde{E}_{v, \mathfrak{p}_{1}^{\infty}})\simeq
Hom(\mathbb{Z}_{p} , \mathbb{Q}_{p}/\mathbb{Z}_{p}) =
\mathbb{Q}_{p}/{\mathbb{Z}}_{p}. \qquad \square\]

Note that there are $r$ primes above $\mathfrak{p}_{2}$ in $F_{n}$
($n=0,1, \ldots)$. It follows from lemma 8 and lemma 9 that
\begin{equation*}
 \mbox{Ker}(\gamma_{n})= \bigoplus_{v \mid p}\mbox{Ker}(\gamma_{n,v})\simeq
\big(\mathbb{Q}_{p}/{\mathbb{Z}}_{p}\big)^{r}.
\end{equation*}
We can now rewrite the exact sequence ({\ref{snake1}) as
\begin{equation}
0 \longrightarrow \mbox{Ker}(\alpha_{n}) \longrightarrow
\mathbb{Q}_{p}/{\mathbb{Z}}_{p} \longrightarrow
\Big(\mathbb{Q}_{p}/{\mathbb{Z}}_{p}\Big)^{r}  \longrightarrow
\mbox{Coker}(\alpha_{n}) \longrightarrow 0. \label{aa}
\end{equation}

The above exact sequence enables us to deduce the following result
about the $\Lambda$-module structure of $X(F_{\infty})$:\\

\noindent {\bf Lemma 10 :} $X(F_{\infty})_{\Gamma_{n}}$ is a free
$\mathbb{Z}_{p}$-module. \\

\noindent {\bf Proof} : Taking the Pontrjagin dual of the exact
sequence ({\ref{aa}}), we obtain
\[ 0 \longrightarrow \widehat{\mbox{Coker}(\alpha_{n}}) \longrightarrow
\mathbb{Z}_{p}^{r} \longrightarrow \ldots \] This tells us that
$\widehat{\mbox{Coker}(\alpha_{n})}$ is a finitely generated free
$\mathbb{Z}_{p}$-module. Taking Pontrjagin dual in the first
column of the commutative diagram (c.d.), we have
\[ 0 \longrightarrow \widehat{\mbox{Coker}(\alpha_{n})}
\longrightarrow X(F_{\infty})_{\Gamma_{n}}\longrightarrow
X(F_{n}). \] By proposition 3, we know that $X(F_{n})$ is a free
$\mathbb{Z}_{p}$-module. As both
$\widehat{\mbox{Coker}}(\alpha_{n})$ and  $X(F_{n})$ have no
$\mathbb{Z}_{p}$-torsion, it is clear that
$X(F_{\infty})_{\Gamma_{n}}$ is a free $\mathbb{Z}_{p}$-module.
\hspace*{1cm} $\square$\\

\noindent{\bf Proof of theorem 1} : \noindent We shall show that
the exact sequence (\ref{aa}) and lemma 10 imply theorem 1. By
considering the $\mathbb{Z}_{p}$-coranks of the terms in the exact
sequence (\ref{aa}), we find that
\begin{equation}
\label{kc}
corank_{{\mathbb{Z}}_{p}}\big(\mbox{Coker}(\alpha_{n})\big)-corank_{{\mathbb{Z}}_{p}}
\big(\mbox{Ker}(\alpha_{n})\big) = r - 1.
\end{equation}
The left vertical map in the commutative diagram (c.d.) implies that
\begin{equation*}
\begin{split}
  &corank_{{\mathbb{Z}}_{p}}\big( Sel_{\mathfrak{p}_{1}^{\infty}}
  (E/F_{\infty})^{\Gamma_{n}}\big)\\
= \;&corank_{{\mathbb{Z}}_{p}}\big(\mbox{Coker}(\alpha_{n})\big)-
corank_{{\mathbb{Z}}_{p}}\big(\mbox{Ker}(\alpha_{n})\big)+
corank_{{\mathbb{Z}}_{p}} \big(Sel_{\mathfrak{p}_{1}^{\infty}}(E/F_{n})\big)\\
= \;&  r - 1 + corank_{{\mathbb{Z}}_{p}}
\big(Sel_{\mathfrak{p}_{1}^{\infty}}(E/F_{n})\big),\;\;[\mbox{by}\;\;(\ref{kc})]
\end{split}
\end{equation*}
i.e.,
\begin{equation}
\label{dis}
 rank_{{\mathbb{Z}}_{p}}\big(X(F_{\infty})\big)_{\Gamma_{n}} =
 \lambda_{n} + r-1.
\end{equation}
By lemma 10, we can conclude that
\begin{equation*}
\big(X(F_{\infty})\big)_{\Gamma_{0}} \simeq
\mathbb{Z}_{p}^{\lambda_{0} + r-1}.
\end{equation*}
In particular, we have
\begin{equation*}
X(F_{\infty})/(p,T) \simeq \big(\mathbb{Z}_{p}/p
\big)^{\lambda_{0}+r-1} = \mbox{a finite module}.
\end{equation*}
 Since $(p,T)$ is the maximal ideal of
 $\mathbb{Z}_{p}[[T]]\simeq \Lambda$, theorem
1 follows from Nakayama's lemma (e.g., see pp. 126 of [La]) for
compact
$\Lambda$-modules. \hspace*{1cm} $\square$\\

\begin{center}7. \underline{\large $\Lambda$-{\it rank of} $X(F_{\infty})$} \end{center}

We have shown in the preceding section that $X(F_{\infty})$ is a
finitely generated $\Lambda$-module. We want to compute its
$\Lambda$-rank and its $\Lambda$-torsion submodule. By structure
theory of $\Lambda$-modules [see (\ref{pan2}) and `General Lemma'
near the end of this section], it will be enough to show that
$(X(F_{\infty}))_{\Gamma_{n}}$ is a free $\mathbb{Z}_{p}$-module
of rank $p^{n}.c$, where $c$ is a constant independent of $n$.
Then, the `General Lemma' would imply that $X(F_{\infty})$ is a
free $\Lambda$-module of rank $c$. Since the $\mathbb{Z}_{p}$-rank
of $(X(F_{\infty}))_{\Gamma_{n}}$ is $(\lambda_{n} + r -1 )$ by
(\ref{dis}), we want to know how the $\lambda_{n}$'s vary with $n$
as we go along the tower of fields
$F_{n}$ over $F_{0}$.\\

\noindent {\bf Lemma 11 :} $\lambda_{n+1} = p\lambda_{n} + (p-1) (r-1) $.\\

\noindent We prove lemma 11 using ideas from [H-M]. Let $G$ be the
Galois group Gal($F_{n+1}/F_{n}$). It is obvious that $G$ is a
cyclic group of order $p$. Formula (3.3) in [H-M] implies that
\begin{equation*}
\begin{split}
& corank_{\mathbb{Z}_{p}}\big( Sel_{\mathfrak{p}_{1}^{\infty}}
(F_{n+1})\big)\\
= & \;p. corank_{{\mathbb{Z}}_{p}}
\big(Sel_{\mathfrak{p}_{1}^{\infty}}(F_{n})\big) +
(p-1)\mbox{ord}_{p} \big(h_{G}(
Sel_{\mathfrak{p}_{1}^{\infty}}(F_{n+1}))\big),
\end{split}
\end{equation*}
where $h_{G}$ denotes the Herbrand quotient. In our notation, the
above formula becomes
\begin{equation}
\label{hm11} \lambda_{n+1} = p.\lambda_{n} + (p-1)\mbox{ord}_{p}
\big(h_{G}( Sel_{\mathfrak{p}_{1}^{\infty}}(F_{n+1}))\big).
\end{equation}
We will now calculate the Herbrand quotient of the Selmer group in
the above expression, since it will determine the explicit
relation between $\lambda_{n+1}$ and $\lambda_{n}$.
 The second exact sequence in the commutative diagram
(c.d.) of section 6 implies that
\begin{equation}
\label{hm1} h_{G}( Sel_{\mathfrak{p}_{1}^{\infty}}(F_{n+1}))=
\frac{h_{G}\big(H^{1}(G(F_{S}/F_{n+1}),
E_{\mathfrak{p}_{1}^{\infty}})\big)} {\prod\limits_{v \mid
p}h_{G}\big(H^{1}(F_{n+1,v},E)_{\mathfrak{p}_{1}^{\infty}}\big)}.
\end{equation}
We shall evaluate the numerator and the denominator in the above
expression with the next three propositions. We shall adopt arguments
of Hachimori and Matsuno who dealt with the cyclotomic situation. The
following lemma simplifies the calculation of the right hand side of
(\ref{hm1}). \\

\noindent{\bf Lemma 12 :} For $i=1,2$, we have
\begin{equation*}
\begin{split}
& (a) H^{i}\big(G,H^{1}(G(F_{S}/F_{n+1}),
E_{\mathfrak{p}_{1}^{\infty}})\big)
=H^{i}(G,E_{\mathfrak{p}_{1}^{\infty}}), \\
 &(b) H^{i}\big(G,H^{1}(F_{n+1,v},E)_{\mathfrak{p}_{1}^{\infty}}\big) =
H^{i}\big(G,E(F_{n+1,v})\big)_{\mathfrak{p}_{1}^{\infty}}.
\end{split}
\end{equation*}

\noindent{\bf Proof : }(a) The Galois group Gal$(F_{S}/F_{n+1})$
has $p$-cohomological dimension at most 2 (see Prop. 8.3.17 in
[N-S-W]). Combining this with corollary 6, we conclude that
$H^{2}(F_{S}/F_{n+1},\;E_{\mathfrak{p}_{1}^{\infty}})$ vanishes
for $i \geq 2$. Then, we have a long exact Hochschild-Serre
spectral sequence
\begin{align*}
\ldots & H^{2}(F_{S}/F_{n},\;E_{\mathfrak{p}_{1}^{\infty}})
\longrightarrow
H^{1}\big(G,\;H^{1}(F_{S}/F_{n+1},\;E_{\mathfrak{p}_{1}^{\infty}})\big)
\longrightarrow H^{3}(G,\;E(F_{n+1})_{\mathfrak{p}_{1}^{\infty}}) \\
 \longrightarrow & H^{3}\big(F_{S}/F_{n},\;E_{\mathfrak{p}_{1}^{\infty}}\big)
 \longrightarrow H^{2}\big(G,\;H^{1}(F_{S}/F_{n+1},\;E_{\mathfrak{p}_{1}^{\infty}})\big)
 \longrightarrow H^{4}\big(G,\;E(F_{n+1})_{\mathfrak{p}_{1}^{\infty}}\big) \\
\longrightarrow &
H^{4}(F_{S}/F_{n},\;E_{\mathfrak{p}_{1}^{\infty}}) \ldots
\end{align*}
As $G$ is a finite cyclic group, we have
\[ H^{i}(G,\;A)= H^{i+2}(G,\;A)\quad \forall i \geq 0,\]
where $A$ is any $G$-module. As
$H^{i}(F_{S}/F_{n},\;E_{\mathfrak{p}_{1}^{\infty}})=0$ for $i\geq
2$, this part of the lemma holds. \\
(b) The Galois group Gal$(\bar{F}_{n+1,v}/F_{n+1,v})$ has strict
cohomological dimension at most 2 (see Prop. 1 and 4, Ch.II in
[Se]). Moreover, $H^{2}(F_{n+1,v},E)$ is trivial because
\begin{equation*}
H^{2}(F_{n+1,v},E) = \lim_{\substack{\longrightarrow\\
\mathbb{Q}_{p}
\subset M \subset F_{n+1,v}\\
[M \;:\; \mathbb{Q}_{p}]<\infty}}H^{2}(M,E),
\end{equation*}
and by Tate local duality (see Ch.II prop. 16 in [Se]),
$H^{2}(M,E)$ vanishes for any finite extension $M$ of
$\mathbb{Q}_{p}$. As in the previous proposition, we have a long
exact Hochschild-Serre spectral sequence and we can conclude that
\begin{equation*}
 H^{i}\big(G,H^{1}(F_{n+1,v},E)\big)_{\mathfrak{p}_{1}^{\infty}}
 =
H^{i+2}\big(G,E(F_{n+1,v})\big)_{\mathfrak{p}_{1}^{\infty}}\;\;\mbox{for}\;\;i=1,2.
\end{equation*}
As $H^{1}(F_{n+1,v},E)$ is a torsion group and $G$ is cyclic, the
above expression reduces to
\[
H^{i}\big(G,H^{1}(F_{n+1,v},E)_{\mathfrak{p}_{1}^{\infty}}\big)
 =H^{i}\big(G,E(F_{n+1,v})\big)_{\mathfrak{p}_{1}^{\infty}} \;\;
\;\;\mbox{for}\;\;i=1,2.\quad \square\]

 \noindent{\bf Proposition 13 :} $h_{G}\big(H^{1}(G(F_{S}/F_{n+1}),
E_{\mathfrak{p}_{1}^{\infty}})\big)=\frac{1}{p}$.\\

\noindent {\bf Proof} : By the first part of lemma 12,
\[h_{G}\big(H^{1}(G(F_{S}/F_{n+1}), E_{\mathfrak{p}_{1}^{\infty}})\big)=
h_{G}(E_{\mathfrak{p}_{1}^{\infty}}).\] Clearly, $G$ acts
trivially on $E_{\mathfrak{p}_{1}^{\infty}}$ as these points are
defined over $F_{n}$. Let $s$ be a generator of $G$ and suppose
$N=\sum_{i=o}^{p-1}s^{i}$. Then
\begin{align*}
H^{2}(G,\;E_{\mathfrak{p}_{1}^{\infty}}) & =
(E_{\mathfrak{p}_{1}^{\infty}})^{G}/N(E_{\mathfrak{p}_{1}^{\infty}}) = 0,\\
H^{1}(G,\;E_{\mathfrak{p}_{1}^{\infty}}) & = \mbox{Ker}(N)/(s-1)
E_{\mathfrak{p}_{1}^{\infty}} = E_{\mathfrak{p}_{1}}.
\end{align*}
Therefore,
\[h_{G}\big(H^{1}(G(F_{S}/F_{n+1}), E_{\mathfrak{p}_{1}^{\infty}})\big)=
h_{G}(E_{\mathfrak{p}_{1}^{\infty}})=\frac{1}{p}. \quad \square\]
\noindent
We calculate the denominator in (\ref{hm1}) by proving the following two propositions.\\

\noindent{\bf Proposition 14 :} $h_{G}\big(H^{1}
(F_{n+1,v},E)_{\mathfrak{p}_{1}^{\infty}}\big)=1 \;\; \forall v
\mid
\mathfrak{p}_{1}$.\\

 \noindent {\bf Proof} : By the second part of lemma 12, we need to calculate the
 ratio of the order of $H^{i}\big(G, E(F_{n+1,v})\big)_{\mathfrak{p}_{1}^{\infty}}$ for
 $i=2,1$.
We consider the following exact sequence of $G$-modules
\begin{equation}
\label{formal} 0 \rightarrow \hat{E}(\mathfrak{m}_{n+1,v})
\rightarrow E(F_{n+1 , v}) \rightarrow
\tilde{E}_{v}(k_{n+1,v})\rightarrow 0,
\end{equation} where $\mathfrak{m}_{n+1,v}$ is the
maximal ideal of $F_{n+1,v}$, and $k_{n+1,v}$ is the residue
field. Taking $G$-cohomology, we have a long exact sequence
\begin{equation}\begin{split}
\label{les} \ldots &\longrightarrow
H^{1}\big(G,\hat{E}(\mathfrak{m}_{n+1,v})\big)_{\mathfrak{p}_{1}^{\infty}}
\longrightarrow
H^{1}\big(G, E(F_{n+1 , v})\big)_{\mathfrak{p}_{1}^{\infty}}\longrightarrow \\
& H^{1}\big(G,\;\tilde{E}_{v}(k_{n
+1,v})\big)_{\mathfrak{p}_{1}^{\infty}} \longrightarrow
H^{2}\big(G
,\hat{E}(\mathfrak{m}_{n+1,v})\big)_{\mathfrak{p}_{1}^{\infty}}
\longrightarrow \ldots.
\end{split}\end{equation}

\noindent For $v|{\mathfrak{p}_{1}}$, $F_{n+1,v}$ is deeply
ramified. By a result of Coates and Greenberg (Theorem 3.1 in
[C-G]), $H^{i}\big(G,\hat{E}(\mathfrak{m}_{n+1,v})\big)=0\;
\forall i \geq 1$. Moreover, $\tilde{E}_{v}(k_{n+1,v})$ is a
torsion group and we can take the
$\mathfrak{p}_{1}^{\infty}$-torsion inside the cohomology in
(\ref{les}). We now have
\begin{equation}
\label{formal2} H^{i}\big(G,
E(F_{n+1,v})\big)_{\mathfrak{p}_{1}^{\infty}}=
H^{i}\big(G,\tilde{E}_{v}(k_{n+1,v})_{\mathfrak{p}_{1}^{\infty}}\big)
\quad \mbox{for} \;\; i=1,2.
\end{equation}
For $v|\mathfrak{p}_{1}$, $k_{n+1,v}$ is the residue field of a
ramified $\mathbb{Z}_{p}$-extension of a finite extension of
$\mathbb{Q}_{p}$, and hence $k_{n+1,v}$ is a finite field. Let us
now consider the $\mathfrak{p}_{1}$-primary part in
(\ref{formal}):
\begin{equation*}
0 \rightarrow
\hat{E}(\mathfrak{m}_{n+1,v})_{\mathfrak{p}_{1}^{\infty}}
\rightarrow E(F_{n+1 , v})_{\mathfrak{p}_{1}^{\infty}} =
\mathbb{Q}_{p}/\mathbb{Z}_{p} \rightarrow
\tilde{E}_{v}(k_{n+1,v})_{\mathfrak{p}_{1}^{\infty}}= \mbox{a
finite module} \rightarrow 0.
\end{equation*}
But $\mathbb{Q}_{p}/\mathbb{Z}_{p}$ has no nontrivial finite
quotient, and
 we deduce that
{$\tilde{E}_{v}(k_{n+1,v})_{\mathfrak{p}_{1}^{\infty}}=0$.
Therefore, $H^{i}\big(G,\tilde{E}_{v}(k_{n
+1,v})_{\mathfrak{p}_{1}^{\infty}}\big)$=0. By lemma 12 (b) and
(\ref{formal2}), we now conclude that
\[H^{i}\big(G,H^{1}(F_{n+1,v},E)_{\mathfrak{p}_{1}^{\infty}}\big) = 0 \;\;\forall
v \mid \mathfrak{p}_{1}\quad \mbox{for}\;\;i=1,2. \] In
particular, the
Herbrand quotient $h_{G}$ is 1. \hspace*{1cm} $\square$\\}

\noindent {\bf Proposition 15 :}
$h_{G}\big(H^{1}(F_{n+1,v},E)_{\mathfrak{p}_{1}^{\infty}}\big)=\frac{1}{p}
\;\;
\forall v \mid \mathfrak{p}_{2}$.\\

\noindent{\bf Proof} : We proceed as in the previous proposition.
However, $\pi$ is an automorphism of $\hat{E}$  for  $v$ not
dividing $\pi$. Therefore,
$H^{i}\big(G,\hat{E}(\mathfrak{m}_{n+1,v})\big)_{\mathfrak{p}_{1}^{\infty}}=0\;\forall\;i\geq
0$. By (\ref{les}),
\begin{equation}
\label{fu2} H^{i}\big(G, E(F_{n+1 ,
v})\big)_{\mathfrak{p}_{1}^{\infty}} =
H^{i}\big(G,\tilde{E}_{v}(k_{n+1,v})\big)_{\mathfrak{p}_{1}^{\infty}}\;\;\forall
i\geq 0.
\end{equation}
As before, we can take the $\mathfrak{p}_{1}^{\infty}$-torsion
inside the cohomology on the right hand side of (\ref{fu2}). Since
the extension $F_{n+1,v}$ is totally ramified over $F_{n,v}$, the
Galois group $G$ acts trivially on $\tilde{E}_{v}(k_{n+1,v})$.
Clearly,
\begin{align*}
\mid H^{1}\big(G,\tilde{E}_{v}(k_{n
+1,v})\big)_{\mathfrak{p}_{1}^{\infty}}\mid & =
\;\mid Hom(G,\mathbb{Q}_{p}/\mathbb{Z}_{p})\mid\; = p \\
\mid H^{2}\big(G,\tilde{E}_{v}(k_{n
+1,v})\big)_{\mathfrak{p}_{1}^{\infty}}\mid & = \;\mid H^{2}(G,
\mathbb{Q}_{p}/\mathbb{Z}_{p}) \mid  = 1.
\end{align*}
From lemma 12 (b) and (\ref{fu2}), it is now obvious that
$h_{G}\big(H^{1}(F_{n+1 , v},E)_{\mathfrak{p}_{1}^{\infty}}\big)=
\frac{1}{p}$.
\hspace*{1cm}$\square$\\

We can now derive the relation between $\lambda_{n+1}$ and
$\lambda_{n}$, as stated in lemma 11. We substitute the values
obtained by the three previous propositions in (\ref{hm1}). We
find that
\begin{equation*}
h_{G}\big(Sel_{\mathfrak{p}_{1}^{\infty}}(F_{n+1})\big) =
\frac{\frac{1}{p}}{(\frac{1}{p})^{r}} = p^{r-1},
\end{equation*}
recalling that $r$ is the number of primes above
$\mathfrak{p}_{2}$ in $F_{n+1}$ for any $n$. Now, it follows from
(\ref{hm11}) that
\begin{equation*}
\lambda_{n+1}= p\lambda_{n}+(p-1)(r-1).
\end{equation*}
This completes the proof of lemma 11.\\

\noindent { \bf Lemma 16 :} $X(F_{\infty})_{\Gamma_{n}}$ is a free
$\mathbb{Z}_{p}$-module of rank $p^{n}(\lambda_{0}+r-1$). \\

\noindent {\bf Proof} : We already saw that
$X(F_{\infty})_{\Gamma_{n}}$ is a free $\mathbb{Z}_{p}$-module
[c.f. lemma 10] of rank $(\lambda_{n}+r-1)$ [c.f. (\ref{dis})]. By
using lemma 11 recursively, we obtain that
\[\lambda_{n}= p^{n}\lambda_{0}+(r-1)(p^{n}-1).\]
Substituting in (\ref{dis}), we find that
\[rank_{{\mathbb{Z}}_{p}}\big(X(F_{\infty})\big)_{\Gamma_{n}} =
p^{n}(\lambda_{0} + r -1). \qquad \square \]

We can now prove Theorem 2 with the following result about the
structure of $\Lambda$-modules (the proof is included for the sake
of completeness):\\

\noindent{\bf General Lemma :} Let $Y$ be a $\Lambda$-module such
that $Y_{\Gamma_{n}}$ is a free $\mathbb{Z}_{p}$-module of rank
$cp^{n}$. Then $Y$ is a free $\Lambda$-module of rank $c$. \\

\noindent {\bf Proof :} Recall that $\Lambda \simeq
\mathbb{Z}_{p}[[T]]$. By structure theory of finitely generated
$\mathbb{Z}_{p}[[T]]$-modules, there is a homomorphism $\psi$ of
$\mathbb{Z}_{p}[[T]]$-modules
\begin{equation}
\label{pan2} \scriptstyle{0 \longrightarrow A \longrightarrow Y
\stackrel{\psi}{\longrightarrow } N = \bigoplus
\mathbb{Z}_{p}[[T]]^{a}\bigoplus \big(\oplus_{i=1}^{s}
\mathbb{Z}_{p}[[T]]/p^{n_{i}}\big) \bigoplus \big(\oplus_{j=1}^{t}
\mathbb{Z}_{p}[[T]]/(f_{j}^{m_{j}})\big) \longrightarrow B
\longrightarrow 0,}
\end{equation}
where $A$ and $B$ are finite. For sufficiently large $n$,
$\Gamma_{n}$ acts trivially on the finite modules $A$ and $B$.
Therefore, $B^{\Gamma_{n}}=B$, $A_{{\Gamma_{n}}}=A$ for $n$
sufficiently large. We can rewrite (\ref{pan2}) as
\begin{equation*}
\begin{split}
0 \longrightarrow A \longrightarrow Y \longrightarrow \mbox{Im}(\psi)
\longrightarrow 0, \\
0 \longrightarrow \mbox{Im}(\psi) \longrightarrow N
\longrightarrow B \longrightarrow 0.
\end{split}\end{equation*}
Therefore, we have exact sequences
\begin{equation}
\label{ex1} \mbox{Im}(\psi)^{\Gamma_{n}}\longrightarrow
A_{\Gamma_{n}} \longrightarrow Y_{\Gamma_{n}} \longrightarrow
\big(\mbox{Im}(\psi)\big)_{\Gamma_{n}}\longrightarrow 0,
\end{equation}
\begin{equation} \label{ex2}
 N^{\Gamma_{n}}\longrightarrow
B^{\Gamma_{n}}\longrightarrow
\big(\mbox{Im}(\psi)\big)_{\Gamma_{n}}\longrightarrow
N_{\Gamma_{n}} \longrightarrow B_{\Gamma_{n}}\longrightarrow 0.
\end{equation}
 By our assumption, it is now clear from (\ref{ex1}) that
 $\big(\mbox{Im}(\psi)\big)_{\Gamma_{n}}$ is a free
 $\mathbb{Z}_{p}$-modules of rank $p^{n}c$. From (\ref{ex2}), we can now deduce
 that $a=c$ and $N$ has no $\mathbb{Z}_{p}[[T]]$-torsion part (note that the order of
 $B_{\Gamma_{n}}$ is bounded independent of $n$). Thus,
 $N= \mathbb{Z}_{p}[[T]]^{c}$. Therefore,
 $N^{\Gamma_{n}}=0$ and $B^{\Gamma_{n}}
 \hookrightarrow \big(\mbox{Im}(\psi)\big)_{\Gamma_{n}}$. Since
 $\big(\mbox{Im}(\psi)\big)_{\Gamma_{n}}$
 does not have any nontrivial finite $\mathbb{Z}_{p}$-submodule,
 $B^{\Gamma_{n}}=0$
 for all $n$. Thus, $B=0$ and $\mbox{Im}(\psi)=N = \mathbb{Z}_{p}[[T]]^{c}$. Now,
 $\mbox{Im}(\psi)^{\Gamma_{n}}=0$, and (\ref{ex1}) implies that
 $A_{\Gamma_{n}}
 \hookrightarrow Y_{\Gamma_{n}}$. But $Y_{\Gamma_{n}}$
 does not have any nontrivial finite $\mathbb{Z}_{p}$-submodule. Thus, $A_{\Gamma_{n}}=0$
 for all $n$. Therefore, $A= 0$. We can now rewrite (\ref{pan2}) as
\[Y \cong \mathbb{Z}_{p}[[T]]^{c}. \qquad \square\]
 \noindent {\bf Proof of Theorem 2 :} Theorem 2 follows directly
from lemma
16 and the `General Lemma' above. \hspace*{1cm}$\square$\\

We can conclude that when $F(E_{p^{n}})$ is abelian over $K$ for
all $n\geq0$, the Pontrjagin dual $X\big(F(E_{{p}^{\infty}})\big)$
of the $\mathfrak{p}_{1}^{\infty}$-Selmer group of $E$ over
$F(E_{{p}^{\infty}})$ is a free $\mathbb{Z}_{p}[[T]]$-module of
rank $\lambda_{0}+r-1$. In particular, it is true when $E$ is
defined over $K$ as the abelian property is implied by theory of
complex multiplication.

\end{document}